\documentclass[12pt]{amsart}
\usepackage{amsmath,amstext,amssymb,amscd}
\usepackage{fullpage}
\usepackage{eucal}

\theoremstyle{plain}
\newtheorem{thm}{Theorem}[section]
\newtheorem{cor}[thm]{Corollary}

\newtheorem{lma}[thm]{Lemma}

\theoremstyle{definition}
\newtheorem{example}[thm]{Example}
\newenvironment{remark}{\par\noindent{\bf Remark.}}{{\vskip 0.095in}}

\def\resultspace{{\vskip 0.075in}}

\def\bbC{\mathbb C}

\def\bbN{\mathbb N}

\def\bbR{\mathbb R}
\def\bbZ{\mathbb Z}
\def\fA{{\mathfrak A}}
\def\sA{{\mathcal A}}

\def\sB{{\mathcal B}}

\def\sC{{\mathcal C}}
\def\fD{{\mathfrak D}}

\def\sD{{\mathcal D}}
\def\sE{{\mathcal E}}
\def\sF{{\mathcal F}}

\def\cG{{\mathcal G}}

\def\sI{{\mathcal I}}

\def\bM{{\mathbf M}}

\def\sP{{\mathcal P}}

\def\sS{{\mathcal S}}

\def\bT{{\mathbf T}}

\def\hull(#1){\mathrm{hull}(#1)}
\def\kker(#1){\mathrm{ker}(#1)}
\def\ker{\mathrm{ker}\thinspace}

\def\FR{\mathrm{R}}

\def\Id{\mathrm{Id\/}}
\def\ordp#1{#1^{\mathrm{ord}}}
\def\nor{N_\fD(\fA)}

\def\ordp#1{#1^{\mathrm{ord}}}

\def\cstar{{$\mathrm{C}^*$}}
\def\Id{\mathrm{Id}}

\def\stbar{{\thinspace |\thinspace }}
\def\dstbar{{\thinspace \bigl|\thinspace }}
\def\indlimit{{\displaystyle \lim_{\longrightarrow} }}

\newcommand{\fe}{\varphi}
\newcommand{\sfe}{\sS(\fe)}
\newcommand{\om}{\omega}

\newcommand{\wh}{\widehat}
\newcommand\ov{\overline}

\newcommand\ur{\underline{r}}
\newcommand\us{\underline{s}}

\def\mapping#1#2#3{{ #1 \colon #2 \rightarrow #3 }}

\begin{document}
\bibliographystyle{amsplain}
\title
       [Algebraic Isomorphisms of Limit Algebras]
       {Algebraic Isomorphisms\\ of Limit Algebras}
\author
       [A. P. Donsig]
       {A. P. Donsig}
\address{Department of Mathematics and Statistics\\
         University of Nebraska at Lincoln\\
         Lincoln, NE, U.S.A  68588-0323}
\email{adonsig@math.unl.edu}
\author
       [T. D. Hudson]
       {T. D. Hudson$^*$}
\thanks{${}^*$Research partially supported by an NSF grant.}
\address{Department of Mathematics \\
         East Carolina University \\
         Greenville, NC, U.S.A  27858--4353}
\email{tdh@math.ecu.edu}
\author
       [E. G. Katsoulis]
       {E. G. Katsoulis}
\address{Department of Mathematics \\
         East Carolina University \\
         Greenville, NC, U.S.A  27858--4353}
\email{katsoulise@mail.ecu.edu}
\thanks{March 17, 1998; revised October 14, 1998.}

\begin{abstract}
We prove that algebraic isomorphisms between limit algebras 
are automatically continuous, and consider the consequences of 
this result.
In particular, we give partial solutions to a conjecture of 
Power~\cite[Notes to Chapter~8]{Power92}, and to an open
problem~\cite[Problem 7.8]{Power92}. 
As a further consequence, we describe epimorphisms between various 
classes of limit algebras. 
%
%
\end{abstract}
\resultspace
\maketitle
\vspace{1cm}

In this paper, we study automatic continuity for limit algebras.
Automatic continuity involves algebraic conditions on a linear 
operator from one Banach algebra into another that guarantee the 
norm continuity of the operator. 
This is a generalization, via the open mapping theorem, of
the uniqueness of norms problem.
Recall that a Banach algebra $\sA$ is said to have a unique 
(Banach algebra) topology if any two complete algebra norms 
on $\sA$ are equivalent, so that the norm topology determined by
a Banach algebra is unique. 
Uniqueness of norms, automatic continuity, and related questions have
played an important and long-standing role in the theory of Banach 
algebras~\cite{Eid40,Sil47,Ric50,Joh67b,Joh67a,DalLoyWil94}. 

Limit algebras, whose theory has grown rapidly in recent years, 
are the nonselfadjoint analogues of UHF and AF \cstar-algebras.
We first prove that algebraic isomorphisms between limit algebras are 
automatically continuous (Theorem~\ref{thm:continuous}). 
This proof uses the ideal theory of limit algebras as well as 
key results from the theory of automatic continuity for Banach algebras. 
Combining this with~\cite[Theorem 8.3]{Power92} verifies 
Power's conjecture that the \cstar-envelope of a limit algebra 
is an invariant for purely algebraic isomorphisms, for limits of
finite dimensional nest algebras, and in particular, for all triangular
limit algebras (Corollary~\ref{cor:power}). 
In~\cite{DonHud96}, the first two authors studied triangular
limit algebras in terms of their lattices of ideals.
By combining automatic continuity with this work, we show that
within the class of algebras generated by their
order preserving normalizers (see below for definitions),
algebraically isomorphic algebras are isometrically isomorphic
(Theorem~\ref{thm:algtoisom}).
This shows that the spectrum, or fundamental relation~\cite{Pow90a},
a topological binary relation which provides coordinates for limit
algebras and is a useful tool in classifications, is a complete algebraic 
isomorphism invariant for this class (Corollary~\ref{cor:alginv}). 
In recent work, the second two authors studied primitivity for limit 
algebras~\cite{HudKatTA}, showing that a variety of limit algebras are 
primitive. 
These results, together with automatic continuity, give  descriptions of 
epimorphisms between various classes of limit algebras, namely
lexicographic algebras (Theorem~\ref{thm:lex}) and
$\bbZ$-analytic algebras (Theorem~\ref{thm:zan}). 

\noindent
\textbf{Preliminaries.}
We briefly recall the framework for studying limit algebras;
see the monograph~\cite{Power92} for details.

Let $\sC$ be an AF \cstar-algebra and $X,Y$ subalgebras 
of $\sC$. 
The {\it normalizer} of $Y$ in $X$ is 
$$
N_Y(X)= 
\bigl\{x\in X \dstbar x\text{ is a partial isometry and }
  xyx^*,x^*yx \in Y \text{ for all } y\in Y \}\,.
$$
A maximal abelian selfadjoint subalgebra (masa) $\sD$ of $\sC$ 
is called a {\it canonical masa} in $\sC$ if there is a nested 
sequence $(C_k)_k$ of finite dimensional \cstar-subalgebras 
of $\sC$ so that $\sC = \overline{\bigcup_k C_k}$, and 
if $D_k = C_k \cap \sD$, then $D_k$ is a masa in $C_k$ 
satisfying $\sD = \overline{\bigcup_k D_k}$ and 
$N_{D_i}(C_i) \subseteq N_{D_{i+1}}(C_{i+1})$.  
A {\it regular canonical subalgebra} is a norm-closed subalgebra of 
an AF \cstar-algebra $\sC$ that contains a canonical masa in $\sC$.
For brevity, we refer to a regular canonical subalgebra as a 
{\it limit algebra}. 
A norm-closed subalgebra $\sA$ of $\sC$ is a limit algebra if 
and only if $\sA$ is the direct limit $\indlimit (A_k,\alpha_k)$ of a 
directed system 
\begin{equation} 
\label{eq:system}
\begin{CD}
A_1 @> \alpha_1 >> A_2 @> \alpha_2 >> A_3 @> \alpha_3 >> A_4 \cdots \,,
\end{CD}
\end{equation}
where for each $k$, 
\begin{itemize}
\item[(i)] $A_k$ is a subalgebra of the finite dimensional 
\cstar-algebra $C_k$ containing a masa $D_k$ of $C_k$,
\item[(ii)] $\alpha_k$ extends to an injective $*$-homomorphism 
from $C_k$ to $C_{k+1}$, and 
\item[(iii)] the extension of $\alpha_k$ maps $N_{D_k}(C_k)$ into 
$N_{D_{k+1}}(C_{k+1})$. 
\end{itemize}
The limit algebra $\sA$ is called {\it triangular AF} (TAF) 
if $\sA \cap \sA^*$ {\it equals} a canonical masa in $\sC$.  
In this case, each algebra $A_k$ in \eqref{eq:system} above 
can be taken to be a subalgebra of some upper triangular matrix 
algebra. 
If $\sA$ is a TAF algebra and $\sA+\sA^*$ is dense in $\sC$, 
then $\sA$ is called {\it strongly maximal triangular}. 
If $\sA$ is a strongly maximal TAF algebra, then we can 
write $\sA= \indlimit (A_k,\alpha_k)$, where each $A_k$ is a 
direct sum of $\bT_n$'s~\cite[Lemma 1.1]{PetPooWag93} ($\bT_n$ 
denotes the $n\times n$ upper triangular matrices with complex 
entries). 
If each $A_k$ is a single $\bT_n$ for each $k$, we call 
the strongly maximal TAF algebra $\sA$ {\it triangular UHF}. 
Note that there are strongly maximal triangular subalgebras of
UHF \cstar-algebras that are not triangular UHF in the above
sense~\cite[Example 2.10]{PetPooWag93}.


\section{Automatic Continuity for Algebraic Isomorphisms of Limit Algebras}

To prove that algebraic isomorphisms of limit algebras are 
continuous, we require Rickart's notion of a separating space. 
Let $\fe :\sA \rightarrow \sB$ be an epimorphism between Banach 
algebras, and define 
$$
\sfe = \bigl\{ b\in\sB \dstbar \text{ there is a sequence }
  (a_n)_n\subseteq \sA \text{ satisfying }
   a_n\rightarrow 0 \text{ and } \fe(a_n) \rightarrow b\bigr\}. 
$$
Obviously $\sfe$ is a closed ideal of $\sB$. 
Moreover, one can easily check that the graph of $\fe$ is 
closed if and only if $\sfe = (0)$.
Thus by the closed graph theorem, $\fe$ is continuous if and only 
if $\sfe = (0)$.  

The following is an adaptation of~\cite[Lemma~2.1]{Sin75}. 

\begin{lma}[Sinclair] 
\label{lma:sinclair}
Let $\fe :\sA\rightarrow \sB$ be an epimorphism of 
Banach algebras, and let $(b_n)_{n=1}^\infty$ be any 
sequence in $\sB$. 
Then there exists $N\in\bbN$ so that for all $n\geq N$, 
$$
\overline{b_1b_2\dots b_n \sfe} = 
   \overline{b_1b_2\dots b_{n+1}\sfe}
$$
and 
$$
\overline{\sfe b_nb_{n-1}\dots b_1} = 
   \overline{\sfe b_{n+1}b_n\dots b_1} \,.
$$
\end{lma}

\begin{proof}
Define $a_n$ by $\fe(a_n)=b_n$. 
The first statement follows from~\cite[Lemma 2.1]{Sin75} with 
$X = \sA$, $Y=\sB$, and $S,T_n,$ and $R_n$ defined as $Sx = \fe(x)$, 
$T_nx = a_nx$, and $R_nx = b_nx$. 
The second statement follows similarly, but with $Sx = \fe(x)$, 
$T_nx = xa_n$, and $R_nx = xb_n$. 
\end{proof}

We also require two technical lemmas about ideals in limit algebras. 
Recall that closed ideals in limit algebras are always 
inductive~\cite{Power92}, that is, if $\sA = \overline{\bigcup_n A_n}$,
and $\sI$ is an ideal, then $\sI = \overline{\bigcup_n (A_n \cap \sI)}$. 

\begin{lma}
\label{lma:1}
If $\sI$ is a nonzero finite dimensional ideal of a limit algebra $\sA$, 
then there exist nonzero projections $p,q$ in $\sA$ so that 
$p\sA q = p\sI q \ne (0)$. 
\end{lma}

\begin{proof} 
Let $\indlimit (A_n,\alpha_n)$ be a presentation for $\sA$, 
where each $A_n$ is a digraph algebra. 
Since $\sI$ is finite dimensional, for some $k$ there are 
finitely many matrix units from $A_k$ which span $\sI$, and 
moreover, for each $n\geq k$, the image under $\alpha_{n,k}$ of each of 
these matrix units is a single matrix unit in $A_n$. 
Let $e$ be a matrix unit in $\sI \cap A_k$, and set 
$p=ee^*$, $q=e^*e$. 
Clearly $p\sI q = \bbC e \ne (0)$. 
Since $e$ has only one restriction in $A_n$, for each $n\geq k$, 
both $p$ and $q$ each have only one restriction in each $A_n$. 
It follows that $p A_n q = p A_k q = p\sI q$ for each 
$n\geq k$, so by inductivity, $p\sA q = p\sI q$.
\end{proof}

\begin{lma} 
\label{lma:2}
If $\sI$ is an infinite dimensional closed ideal in a limit algebra $\sA$, 
then there exists a sequence $(p_n)_{n=1}^\infty$ of mutually orthogonal 
projections in $\sA$ so that either $p_n\sI \ne (0)$ for all $n$ or 
$\sI p_n \ne (0)$ for all $n$.  
\end{lma}
 
\begin{proof} 
Let $\indlimit (A_n,\alpha_n)$ be a presentation for $\sA$. 
By relabeling the presentation, if necessary, we can assume 
$\sI\cap A_1 \ne (0)$. 
We say $e \in N_{\sD}(\sA)$ has infinitely many restrictions if there
is an infinite sequence of pairwise orthogonal projections, $(q_n)$,
in $\sD$ so that $q_n e \ne 0$ for all $n$. 
Note that, if there is a matrix unit $e$ in $\sI$ with infinitely 
many restrictions, then we can choose the sequence $(p_n)$ from the 
orthogonal range projections of restrictions of $e$. 
%
Let $\Omega$ be the collection of all final projections of 
matrix units from $\sI$ and $\Lambda$ the collection of 
initial projections. 
At least one of $\Omega,\Lambda$ must be infinite, for else $\sI$ 
would be finite dimensional.  

Assume $\Omega$ is infinite, and let $\sP_n$ be the set of 
distinct final projections of all matrix units in $\sI\cap A_n$.
If some $p \in \sP_1$ had infinitely many restrictions, then a matrix unit 
in $\sI \cap A_1$ would have infinitely many restrictions and the conclusion
would follow, as noted above.
So we may assume that for every $n$, each $p \in \sP_n$ has only finitely
many restrictions.
Pick some element $p_1 \in \sP_1$.
As $p_1 \in \sP_1 \subseteq \Omega$, $p_1 \sI \ne (0) $.
Since $\sP_1$ is finite and $\Omega$ is infinite, there is some $n > 1$
so that $\sP_n$ contains an element, call it $p_2$, which is not a 
restriction of some $p \in \sP_1$.
Since it is not a restriction of any element of $\sP_1$, $p_2$ may be
chosen orthogonal to each element of $\sP_1$, and in particular, to $p_1$.
Again, $p_2 \sI \ne (0)$.
Next, there is some $m > n$ so that $\sP_m$ contains an element,
call it $p_3$, which is not a restriction of some $p \in \sP_2$.
Continuing in this way, we can construct the required sequence 
$(p_n)_{n=1}^\infty$ with $p_n \sI \ne (0)$.

With the obvious modifications, the same argument works if
$\Lambda$ is infinite.
\end{proof} 

We can now prove our main result. 

\begin{thm}
\label{thm:continuous}
If $\fe :\sA \rightarrow \sB$ is an algebraic isomorphism between 
limit algebras $\sA$ and $\sB$, then $\fe$ is continuous.
\end{thm}

\begin{proof}
It suffices to prove that $\sfe=(0)$. 
Assume that $\sfe$ is nonzero and finite dimensional. 
By Lemma~\ref{lma:1}, there are nonzero projections  
$p,q\in \sB$ so that $p\sB q = p\sfe q \ne (0)$. 
In particular, $p\sB q$ is finite dimensional. 
Let $e,f\in \sA$ satisfy $\fe(e) = p$ and $\fe (f) = q$. 
Clearly $e,f$ are idempotents in $\sA$, and 
$e\sA f$ is finite dimensional since $p\sB q$ is.  

Now let $b\in \sfe$ so that $pbq\ne 0$. 
Then there is a sequence $(a_n)_n$ in $\sA$ so that 
$a_n\rightarrow 0$ and $\fe(a_n) \rightarrow b$. 
Hence 
$$
\fe(ea_nf) = p\fe(a_n)q \longrightarrow pbq \ne 0\,.
$$
But since $ea_nf \rightarrow 0$, this shows that 
the map $\fe|_{e\sA f}$ is discontinuous. 
But dim$\/(e\sA f) < \infty$, a contradiction. 

Hence if $\sfe$ is nonzero, it must be infinite dimensional.  
Applying Lemma~\ref{lma:2} to the ideal $\sfe$, we obtain 
a sequence $(p_n)_{n=1}^\infty$ of mutually orthogonal projections 
in $\sB$ so that either $p_n\sfe \ne (0)$ for each $n$ or 
$\sfe p_n\ne (0)$ for each $n$. 
Set
$$
b_n = \sum_{k=n}^\infty \frac{p_k}{2^k} \,,
$$
so that 
$$
b_1b_2\dots b_n = b_n b_{n-1} \dots b_1 =
	\sum_{k=n}^\infty \frac{p_k}{2^{nk}} \,.
$$
If $p_n\sfe \ne (0)$ for each $n$, then 
$$
\ov{b_1b_2\dots b_n\sfe} \supsetneq 
\ov{b_1b_2\dots b_{n+1}\sfe} \,.
$$
If $\sfe p_n$ is nonzero for each $n$, then 
$$
\ov{\sfe b_nb_{n-1}\dots b_1} \supsetneq 
\ov{\sfe b_{n+1}b_n\dots b_1} \,.
$$
In either case, we have contradicted Theorem~\ref{lma:sinclair}.
It follows that if $\fe$ is an algebraic isomorphism, then 
$\sfe = (0)$, i.e., $\fe$ is continuous.
\end{proof}

\begin{remark}
A similar result for triangular subalgebras of groupoid \cstar-algebras was 
previously announced and circulated in preprint form~\cite[Theorem~6.4]{QiuPR}.
Unfortunately, the proof given there is incomplete.
In any case, our techniques are quite different and 
Theorem~\ref{thm:continuous} applies to non-triangular algebras, albeit
subalgebras of AF \cstar-algebras.
\end{remark}

The above proof, with only trivial changes, also establishes 
the following result.

\begin{cor}
\label{cor:tuhf} 
If $\sA$ is a Banach algebra and $\sB$ a limit algebra with 
no finite dimensional ideals, then any epimorphism from 
$\sA$ onto $\sB$ is continuous.
\end{cor}

Power has shown that the $C^*$-envelope of certain limit algebras
is a Banach algebra invariant~\cite[Theorem 8.3]{Power92}.
Note that Theorem~8.3 is not valid in the generality stated, as
it depends on Lemma~8.1, which is false~\cite{Perrata}; 
thus, Theorem~8.3 must be restricted to inductive limits of 
finite dimensional nest algebras.
However, this includes all triangular limit algebras, so the
theorem still applies to a wide range of algebras.

Power has suggested that the $C^*$-envelope is an invariant for 
purely algebraic isomorphisms of limit 
algebras~\cite[Chapter~8 Notes]{Power92}. 
Power's work and Theorem~\ref{thm:continuous} directly imply
this conjecture for a variety of limit algebras

\begin{cor} 
\label{cor:power} 
If $\sA$ and $\sB$ are algebraically isomorphic limit algebras, both
inductive limits of finite dimensional nest algebras, then 
$C^*(\sA)$ and $C^*(\sB)$ are isomorphic as \cstar-algebras.
\end{cor}

Moreover, Theorem~\ref{thm:continuous} implies that to show the
$C^*$-envelope is an algebraic invariant for all limit algebras,
it is enough to show that it is a Banach algebra invariant.


\section{Algebraic Isomorphisms between Order Preserving TAF Algebras}

Power has conjectured that algebraic isomorphism may be equivalent to
isometric isomorphism for subalgebras of AF \cstar-algebras containing
a regular canonical masa~\cite[p.\ 95]{Pow92c}.
In this section we combine Theorem~\ref{thm:continuous} with results
from~\cite{DonHud96} to prove this conjecture for strongly maximal 
TAF algebras generated by their order preserving normalizers.
That is, two algebras in this family are algebraically 
isomorphic if and only if they are isometrically isomorphic. 

This family includes a range of fundamental examples,
such as the standard, refinement, and alternation limit
algebras~\cite{Bak90,HopPow92,PetPooWag90,Poo92a},
the lexicographic algebras~\cite{PetPoo97,Pow95b,Pow96b}, and
the $\bbZ$-analytic algebras~\cite{PetPooWag93,PooWag93}, as well as
non-analytic algebras~\cite{Hud94,SolVen92}.
We include the necessary definitions below, but refer the reader 
to~\cite{DonHop95,DonHud96} for further details on order 
preservation and TAF algebras generated by their order preserving 
normalizers.

Fix a TAF algebra $\sA$, and let $\sP(\sA)$ denote the collection of 
all diagonal projections of $\sA$.
The {\it diagonal order} on $\sP(\sA)$, denoted ``$\preceq$'', is a
partial order given by
$$
e\preceq f \Longleftrightarrow \text{ there exists }w\in N_\sD(\sA)
   \text{ with }ww^*=e, w^*w=f \,.
$$
Each $w\in N_\sD(\sA)$ induces a partial homeomorphism on $\sP(\sA)$,
which has domain $\{x\in\sP(\sA)\stbar x\leq ww^*\}$ and range
$\{x\in\sP(\sA)\stbar x\leq w^*w\}$, given by $x\longmapsto w^*xw$.
We say that $w$ is {\it order preserving} if this map preserves the
diagonal order restricted to its domain and range.
Define the {\it order preserving normalizer} of $\sA$ to be
$$
N_\sD^{\mathrm{ord}}(\sA) =
        \{w\in N_\sD(\sA)\stbar w\text{ is order preserving}\}\,,
$$
and let $\ordp\sA$ be the subalgebra of $\sA$ generated by
$N^{\mathrm{ord}}_\sD(\sA)$.
Note that a sum of elements of $N_\sD^{\mathrm{ord}}(\sA)$ may 
not be in $N_\sD^{\mathrm{ord}}(\sA)$, even if the sum is in $\nor$.

An embedding $\alpha : A\rightarrow B$ between triangular algebras is called
{\it locally order preserving} if $\alpha(e)$ is an order preserving element
of $B$ for every matrix unit $e \in A$.
For example, standard embeddings, refinement embeddings, embeddings
induced by ordered Bratteli diagrams (see~\cite[Chapter~6]{Power92})
and the block standard embeddings of~\cite{Hud94}
are all locally order preserving.
This is easy to verify using the following characterization.

\begin{thm}[\protect{\cite[Theorem~18]{DonHop95}}]
\label{thm:DH}
A TAF algebra $\sA$ is generated by its order preserving normalizer,
if and only if, it has a presentation $\indlimit(A_i, \alpha_i)$ with
$\alpha_j  \circ \cdots \circ \alpha_i$ locally order preserving 
for all $i\le j$.
\end{thm}

The spectrum, or fundamental relation, of a limit algebra was
introduced by Power~\cite{Pow90b} and is a complete isometric 
isomorphism invariant for TAF algebras.
For a limit algebra $\sA$, let $H$ denote the maximal ideal
space of the canonical masa $\sD \subseteq \sA$.
The {\it spectrum} of $\sA$ is the topological binary relation
$\FR(\sA)$ on $H$ such that $(x,y)\in\FR(\sA)$ if and only if
there is an element $v$ of $N_\sD(\sA)$ so that $y(d)=x(vdv^*)$
for all $d$ in $\sD$.
In this case we say that $v$ {\it relates} $x$ and $y$.
If $e\in N_\sD(\sA)$, define $G(e)$ to be the graph of $e$ in $\FR(\sA)$,
that is,
$$ G(e) = \left\{ (x,y) \in \FR(\sA) \left| \hbox{ $e$ relates $x$ and $y$}
   \right. \right\} \,.
$$
Consider the topology on $\FR(\sA)$ that has
$\{ G(e) \stbar e\in N_\sD(\sA) \}$ as a base of open sets.
The sets $G(e)$ are then clopen and compact, and so $\FR(\sA)$
is a totally disconnected, locally compact Hausdorff space.
Furthermore, the topology on $\FR(\sA)$ is completely determined by those
$G(e)$ for which $e$ is a matrix unit~(see~\cite{Power92} for details).


We will need the following result.

\begin{thm}[\protect{\cite[Corollary 16]{DonHud96}}]
\label{thm:induced}
Let $\sA$ and $\sB$ be strongly maximal TAF algebras generated by their
order preserving normalizers, and assume that the ideal lattices 
of $\sA$ and $\sB$ are isomorphic.
\begin{itemize}
\item[(i)] There is a bijective isometry $\mapping \eta \sA \sB$ and 
closed triangular subalgebras $\sE$ and $\sF$ of $\sA$ with
$\sE + \sF = \sA$ and $\sE\cap \sF$ the diagonal of $\sA$,
so that $\eta$ is an algebra isomorphism on $\sE$ and
an anti-isomorphism on $\sF$.

\item[(ii)] If $\mapping {\phi} {\Id(\sA)} {\Id(\sB)}$ is a lattice 
isomorphism, then $\phi$ is induced by a bijective isometry as 
in {\em (i)}.
\end{itemize}
\end{thm}

The key provided by automatic continuity is that an algebraic 
isomorphism between limit algebras induces a lattice isomorphism 
between their lattices of closed ideals. 
Theorem~\ref{thm:induced} (ii) then yields a bijective isometry 
between the algebras, which is then deduced to be an isomorphism. 

\begin{thm}
\label{thm:algtoisom}
Suppose $\sA$ and $\sB$ are strongly maximal TAF algebras generated
by their order preserving normalizers.
Then $\sA$ and $\sB$ are algebraically isomorphic if and only if they are
isometrically isomorphic.
\end{thm}

\begin{proof}
One direction is trivial; to prove the other, suppose that 
$\mapping \phi \sA \sB$ is an algebraic isomorphism.
We prove the theorem for $\sA$ and $\sB$ unital; the non-unital 
case then follows by adjoining units and applying the unital case.

Let $\sC = \sA \cap \sA^*$.
Theorem~\ref{thm:continuous} shows that an algebra isomorphism is 
automatically continuous, and so $\phi$ maps closed ideals to closed ideals.
Thus it induces a lattice isomorphism which, by 
Theorem~\ref{thm:induced} (ii) is induced by a bijective isometry 
$\mapping \eta \sA \sB$ with the properties described in 
Theorem~\ref{thm:induced} (i).  
We will show that $\sF = \sC$ and so $\sE=\sA$, proving the theorem.

Let $\sD = \sB \cap \sB^*$.
Observe that $\phi$ and $\eta$ induce the same isomorphism from the maximal 
ideals of $\sC$ to the maximal ideals of $\sD$; we claim that
$\eta |_\sC$ and $\phi |_\sC$ are equal.
To see this, first note that these maximal ideal spaces are 
compact Hausdorff spaces (as $\sA$ and $\sB$ are unital). 
Since $\sC$ and $\sD$ are isomorphic to continuous functions on their maximal 
ideal spaces, it follows that $\eta |_\sC$ and $\phi |_\sC$ are both given by
$f \mapsto f \circ \Phi^{-1}$, where $\Phi$ is the common induced map
between the maximal ideal spaces (see \cite[Theorem~3.4.3]{KadisonRin83}).

Suppose $x \in N_\sC(\sA)$ and $x \in \sF\backslash \sC$. 
Notice that since $\eta$ is an isometric algebra anti-isomorphism on $\sF$, 
\cite[Proposition 7.1]{Power92} implies that $\eta(x) \in N_\sD(\sB)$. 
We now claim that the final projection of $\eta(x)$ is $\eta(x^*x)$. 
Indeed, since $x(x^*x) = x$, we have $\eta(x^*x)\eta(x) = \eta(x)$ 
and so $\eta(x)\eta^*(x) \subseteq \eta(x^*x)$. 
On the other hand, for any projection $p\in\sD$ with $p\eta(x) = 0$, 
$x\eta^{-1}(p) = 0$ and so $\eta^{-1}(p) \subseteq 1-x^*x$. 
This implies that $p\subseteq 1-\eta(x^*x)$ or 
$\eta(x^*x) \subseteq 1-p$. 
The reverse inclusion now follows with $p=1-\eta(x)\eta^*(x)$. 
The claim now shows that every point in the graph 
$G\left(\eta(x^*x)\right)$ appears 
as the first element of some ordered pair in $G(\eta(x))$. 
Similarly, the initial projection of $\eta(x)$ is $\eta(xx^*)$ and so 
every point in the graph $G\left(\eta(xx^*)\right)$ appears as the 
second element of some ordered pair in $G(\eta(x))$. 

On the other hand, the element $\phi(x)$ need not be in $N_\sD(\sB)$. 
However, the relation
$$
\phi(xx^*)\phi(x) \phi(x^*x) = \phi(x) \ne 0
$$
implies that $\phi(xx^*)\sB \phi(x^*x) \ne 0$ and so there exists 
non-diagonal $v\in N_\sD(\sB)$ so that 
$G(v^*v) \subseteq G\left(\phi(x^*x)\right)$ and 
$G(vv^*) \subseteq G\left(\phi(xx^*)\right)$. 
Notice that every point in the graph $G(vv^*)$ of the final projection 
of $v$ appears as the first element of some ordered pair in $G(v)$ 
and similarly every point in $G(v^*v)$ appears as the second element of 
some pair in $G(v)$. 

We claim that there is some maximal ideal $z$ of $\sD$ so that 
$(z,z)\in G(v^*v)$ and there is no maximal ideal $y\ne z$ with 
$(y,z) \in \FR(\sB)$ and $(y,y) \in G(v^*v)$. 
Accepting this claim for the moment we obtain a contradiction. 
Since $(z,z) \in G(v^*v)$ the previous paragraph shows that there 
is some maximal ideal $y\ne z$ so that $(y,z) \in G(v)$. 
Since 
$$
(y,y) \in G(vv^*)\subseteq G\left(\phi(xx^*)\right) = G\left(\eta(xx^*)\right)\,,
$$ 
there is a maximal ideal $z^\prime$ so that 
$(z^\prime,y)\in G\left(\eta(x)\right)$ and 
$(z^\prime,z^\prime) \in G\left(\eta(x^*x)\right)$. 
If $z^\prime = z$, then we have 
$(z,y),(y,z)\in \FR(\sB)$ and $z\ne y$, contradicting the anti-symmetry 
of $\FR(\sB)$. 
If $z^\prime \ne z$, then by transitivity $(z^\prime,z)\in \FR(\sB)$, 
contradicting the properties of $z$. 
Hence no such $x$ in $N_\sC(\sA)\cap (\sF\backslash \sC)$ exists and
so $\sF = \sC$.

It remains only to show that we can construct a maximal ideal with the 
desired properties.
Fix a presentation of $\sB$, $\indlimit(B_n,\beta_n)$.
As $v^*v$ is a projection, there is some $n$ so that $v^*v \in B_n$.
Hence for each $m \ge n$, $v^*v$ is a sum of diagonal matrix units 
in $B_m$; let $X_m$ be the set of such diagonal matrix units.
Let $s_m$ be the smallest element of $X_m$ with respect to the diagonal
order; since each $B_m$ is a direct sum of $T_n$'s, it follows that
$s_{m+1}$ is a subordinate of $s_m$ (otherwise, the predecessor of $s_{m+1}$
in $X_m$ is smaller than $s_m$).
By~\cite[Lemma 1]{DonHud96}, there is a bijection between maximal 
ideals in $G(v^*v)$ and sequences $(x_m)_{m \ge n}$ where 
each $x_m \in X_m$ and $x_{m+1}$ is a subordinate of $x_m$.
Since the $s_m$ are minimal in the diagonal order, it follows that the 
maximal ideal in $G(v^*v)$ corresponding to $(s_m)$ satisfies the
conditions of the claim. 
\end{proof}

For triangular limit algebras, the spectrum is an isometric isomorphism 
invariant~\cite[Corollary~7.7]{Power92}, and it has been asked
if this extends to algebraic isomorphism~\cite[Problem 7.8]{Power92}.
Theorem~\ref{thm:algtoisom} combined with the cited result 
yields a partial solution to this problem.  

\begin{cor}
\label{cor:alginv}
The spectrum or fundamental relation $\FR(\sA)$ is a complete 
algebraic isomorphism invariant for strongly maximal TAF algebras 
$\sA$ generated by their order preserving normalizers.
\end{cor}

In addition to this result, automatic continuity can be used 
to provide further evidence for a positive solution to the problem.
That is, we can show that certain not isometrically isomorphic TAF
algebras are not algebraically isomorphic, even though some of 
these algebras are not generated by their order preserving normalizers.

\begin{example}
Let $A_n = \bT_{2^n}$ and define the multiplicity 2 
{\it refinement embedding} $\rho_k: A_k \rightarrow A_{k+1}$ by 
$$
\rho_k([a_{ij}]) = \left[a_{ij}I_2\right]\,,
$$ 
where $I_2$ is the $2\times 2$ identity matrix. 
Now define the {\it elementary twist embedding} 
$\tau_k: A_k\rightarrow A_{k+1}$ by 
$$
\tau_k= (\text{Ad}\,U_{2^{n+1}}) \circ \rho_k\,,
$$
where $U_{2^{n+1}}$ is the permutation unitary in $\bM_{2^{n+1}}$ 
which interchanges the last two minimal diagonal projections of 
$A_{n+1}$. 
Let $\sA = \indlimit(A_n,\rho_n)$ and $\sB = \indlimit(A_n,\tau_n)$.  
The TAF algebras $\sA$ and $\sB$ are called the $2^\infty$ 
{\it refinement } and {\it refinement with twist} algebras, 
respectively. 

One can tell that there is no isometric isomorphism between $\sA$ 
and $\sB$ by examining their spectrums~\cite{Pow90b}. 
But a more delicate question is whether there is an algebraic 
isomorphism between $\sA$ and $\sB$.
If there were such an isomorphism $\phi:\sA\rightarrow\sB$, 
Theorem~\ref{thm:continuous} shows that $\phi$ would be continuous. 
Hence, since $\phi$ then maps closed ideals of $\sA$ to closed 
ideals of $\sB$, $\phi$ induces a lattice isomorphism between the 
lattices of closed ideals of $\sA$ and $\sB$. 
But the sentence before Corollary 28 in~\cite{DonHud96} shows there 
exists no such lattice isomorphism between the lattices of ideals 
of the $2^\infty$ refinement and the refinement with twist, and hence 
there is no algebraic isomorphism between $\sA$ and $\sB$.  
\end{example}
\vspace{0.1cm}

Using an argument similar to the proof of~\cite[Corollary 28]{DonHud96}, 
one can show that there are no algebraic isomorphisms between pairs of 
more general ``twist'' type algebras, as considered 
in~\cite[Section 5]{DonHud96}. 
\vspace{0.5cm}

\section{Epimorphisms}
The class of TAF algebras generated by their order preserving normalizers 
includes all lexicographic algebras and all $\bbZ$-analytic algebras. 
We can describe many of the epimorphisms for these two classes of algebras, 
by using the results of the second two authors on primitivity~\cite{HudKatTA}. 

We start with the lexicographic algebras and recall their construction
from~\cite{Pow95b}.
Let $(\Omega,\leq)$ be a countable linear ordering, i.e., let
$\Omega$ be a countable set and $\leq$ a linear order on $\Omega$.
Let $\omega : \bbN \rightarrow \Omega$ be an enumeration of
$\Omega$ and let $\nu : \Omega \rightarrow \{2,3,4,\dots\}$ be a
multiplicity function;
set $\omega_k = \omega(k)$ and $\nu_k = \nu(\omega_k)$.
Also let $F_n = \{ \omega_1,\omega_2,\dots,\omega_n\}$, for $n\in\bbN$.
View $\bT_{\nu_1\nu_2\dots\nu_n}$ as the subalgebra $\sA_n$ of
$\bM_{\nu_1\nu_2\dots\nu_n}$ which is spanned by the matrix
units
$$
e_{ij} = e_{(i_1,i_2,\dots,i_n),(j_1,j_2,\dots,j_n)}\,,
$$
where $i\leq j$ if and only if
%
%
there is some $k_0 \in\{ 1,2,\dots n\}$, depending on $i$ and $j$,
such that $i_{k_0} < j_{k_0}$ and $i_k = j_k$ for all
$\omega_k <\omega_{k_0}$
(notice that the indexing here is different than that of
Section 2).
%
%
Viewing $\sA_{n+1} = \bT_{\nu_1\nu_2\dots\nu_{n+1}}$ in
a similar manner, we define the map $\psi_n : \sA_n \rightarrow
\sA_{n+1}$ to be the linear extension of the map which takes
$$
e_{ij} \mapsto
\sum_{\xi=1}^{\nu_{n+1}}
e_{(i_1,i_2,\dots,i_n,\xi),(j_1,j_2,\dots,j_n,\xi)} \,.
$$
It is clear that each $\psi_n$ is a $*$-extendible, regular
embedding, and so the direct limit defines a triangular
UHF algebra $\indlimit(\sA_n,\psi_n)$, 
which we denote $A(\Omega,\nu)$.
A standard argument using intertwining diagrams shows that 
a different enumeration of $\Omega$ will produce an 
isometrically isomorphic algebra, so $A(\Omega,\nu)$ is well-defined.

The algebras $A(\Omega,\nu)$, the lexicographic algebras,
were studied in~\cite{Pow95b,Pow96b} and were extended to 
subalgebras of direct sums of $\bM_n$'s in~\cite{PetPoo97}.
The familiar standard, refinement, and alternation
triangular UHF algebras are all examples of lexicographic
algebras.
Indeed, $\Omega = \bbZ_-,\bbZ_+$, and $\bbZ$ yields the
standard, refinement, and alternation algebras, respectively.

Assume that $\Omega = \{ \om_k \stbar k \in \Omega_0\}$ is ordered 
like a subset of the integers, i.e., 
$\Omega_0 \subseteq \bbZ\setminus\{0\}$ and $\om_k\leq \om_l$ 
if $k\leq l$. 
Let $\nu :\Omega\rightarrow \{2,3,4,\dots\}$ be as above. 
We write 
$$
r_k = \nu(\om_k) \quad \text{ and }\quad s_k = \nu(\om_{-k}), \,\,\,
k=1,2,3,\dots
$$
with the understanding that if $k\notin\Omega_0$, then $r_k = 1$, 
and, similarly, if $-k\notin\Omega_0$, then $s_k = 1$. 
Associate with $\mu$ the pair $(\ur,\us)$ of (possibly) generalized 
integers 
$$
\ur = r_1r_2r_3\dots \quad\text{ and }\quad 
\us = s_1s_2s_3\dots 
$$
Define an equivalence relation $\sim$ on pairs of generalized 
(or finite) integers by $(\ur,\us) \sim(\ur^\prime,\us^\prime)$ 
if and only if $\ur\us = \ur^\prime\us^\prime$ and there exist 
coprime natural numbers $a,b$ so that $b\ur= a\ur^\prime$ 
and $a\us = b\us^\prime$. 

Assume now that $\Omega$ is an arbitrary countable linear 
ordering. 
Define an equivalence relation $\approx$ on $\Omega$ 
such that $\om \approx\om^\prime$ if the order intervals 
$[\om,\om^\prime]$ and $[\om^\prime,\om]$ are finite. 
Then the set $\Omega\bigl/\bigr.\approx$ of equivalence 
classes is linearly ordered and each equivalence class  
$\langle\om\rangle$ is itself a linearly ordered set 
isomorphic to a subset of the integers. 
Then to each equivalence class $\langle\om\rangle$ we can 
associate a pair of generalized integers
$p_\nu \bigl(\langle\om\rangle\bigr) = (\ur,\us)$.

Power's classification of lexicographic algebras~\cite[Theorems~4 
and~6]{Pow96b}, after invoking Theorem~\ref{thm:algtoisom}, becomes:

\begin{thm}
\label{thm:Power}
Let $\Omega,\Omega^\prime$ be countable linear orderings with 
maps $\nu:\omega\rightarrow \{2,3,4,\dots\}$ and 
$\nu^\prime :\omega^\prime\rightarrow \{2,3,4,\dots\}$. 

Then $A(\Omega,\nu)$ and $A(\Omega^\prime,\nu^\prime)$ 
are isomorphic as complex algebras if and only if there is an 
order bijection $s:\Omega/\approx \rightarrow  \Omega^\prime/\approx$ 
so that 
$p_{\nu^\prime}\bigl(s\bigl(\langle\om\rangle\bigr)\bigr) 
    = p_\nu\bigl(\langle\om\rangle\bigr)$.
\end{thm}

By~\cite{Pow95b}, the semisimple lexicographic algebras are precisely
those $A(\Omega',\nu')$ where $\Omega'$ has no minimal element.
Further, by~\cite[Proposition 3.1]{HudKatTA}, these algebras
all primitive, so this class is precisely the primitive lexicographic
algebras.

\begin{thm}
\label{thm:lex}
There is an epimorphism from $A(\Omega,\nu)$ onto 
$B$, a primitive Banach algebra, if and only if there is an order 
interval decomposition $\Omega= \wh{\Omega}\oplus\wh{\Omega}^{\text{c}}$,
where $\Omega^{\text{c}}$ contains no minimal element, 
so that 
$A(\wh{\Omega}^{\text{c}},\nu |_{\wh{\Omega}^{\text{c}}})$ is isomorphic
to $B$, as complex algebras.

In particular, if $B=A(\Omega',\nu')$ where, necessarily,
$\Omega^\prime$ has no minimal element, then there is an order bijection 
$s:\wh{\Omega}^{\text{c}}\bigl/\bigr.\approx\rightarrow
   \Omega^\prime\bigl/\bigr.\approx $
so that $p_{\nu^\prime}\bigl(s\bigl(\langle\om\rangle\bigr)\bigr) 
   = p_\nu\bigl(\langle\om\rangle\bigr)$, 
where $\om \in \wh{\Omega}^{\text{c}}$. 
\end{thm}

\begin{proof} 
If $\fe: A(\Omega,\nu)\rightarrow B$ 
is an epimorphism, then $\ker\fe\subseteq A(\Omega,\nu)$ is a 
primitive ideal. 
Then by~\cite[Theorem~3.6]{HudKatTA}, there exists an order 
interval decomposition $\Omega= \wh{\Omega}\oplus\wh{\Omega}^{\text{c}}$, 
where $\wh{\Omega}^{\text{c}}$ has no minimal element, so
that, as a complex algebra, $B$ is isomorphic to 
$$ A(\Omega,\nu)\bigl/\bigr.\ker\fe \approx 
	A(\wh{\Omega}^{\text{c}},\nu  |_{\wh{\Omega}^{\text{c}}})\,.$$
If $B = A(\Omega^\prime,\nu^\prime)$, then the conclusion 
follows from Theorem~\ref{thm:Power}.
\end{proof} 

Theorem~\ref{thm:lex} is false if $B$ is not primitive.
Indeed, it is easy to see that there is an epimorphism from the 
$2^\infty$ refinement algebra onto the $3\cdot 2^\infty$ refinement.

Finally, we investigate the epimorphisms of $\bbZ$-analytic algebras. 
The most convenient definition of $\bbZ$-analyticity involves a
cocycle on the spectrum, so we define cocycles precisely.
%

Recall the definition of spectrum from just before Theorem~\ref{thm:induced}. 
Let $\sB$ be the UHF \cstar-algebra corresponding to $\sA$, and
let $\cG$ be the equivalence relation generated by $\FR(\sA)$.
A continuous function $\mapping c \cG \bbR$ is called a {\it cocycle}
if $c(x,y) + c(y,z) = c(x,z)$ for all points $(x,y), (y,z)\in\cG$.
The algebra $\sA$ is {\it analytic} if there is a cocycle $c$ on
$\cG$ so that $c^{-1}([0,\infty)) = \FR(A)$.
The standard, refinement, and alternation limit algebras are
analytic~\cite[Examples 6.1--6.3]{Ven91};
for strongly maximal non-analytic algebras,
see~\cite[Section~4]{MasMuhSol94}, \cite[Proposition 10.18]{Power92},
or~\cite[Example 3.2]{SolVen92}.
Algebras that are analytic by an integer-valued cocycle are called
{\it $\bbZ$-analytic}~\cite{PetPooWag93,PooWag93}.
Although standard embedding algebras are the generic examples
of $\bbZ$-analytic algebras, the paper~\cite{PetPooWag93} gives an example,
due to Donsig and Hopenwasser, of a semisimple, $\bbZ$-analytic algebra 
which is not built up out of standard embeddings.

\begin{thm}
\label{thm:zan} 
If $\fe :\sA\rightarrow\sB$ is an epimorphism from a $\bbZ$-analytic 
triangular UHF algebra $\sA$ onto a triangular UHF algebra $\sB$,  
then $\fe$ is an isomorphism. 
\end{thm} 

\begin{proof} 
Since $\sB$ is triangular UHF, the zero ideal of $\sB$ is meet 
irreducible. 
By automatic continuity, the closed ideals of $\sA$ containing $\ker \fe$ 
form a sublattice of the closed ideals of $\sB$.
Hence $\ker\fe$ is meet irreducible in $\sA$. 
But $\sA$ is $\bbZ$-analytic, so 
by~\cite[Proposition 4.6]{DonHopHudLamSolTA}, $\ker\fe$ is either 
zero or it must have finite codimension in $\fA$. 
But if $\ker\fe$ had finite codimension, since 
$\sA\left/\right.\ker\fe\simeq\sB$, we would have that 
$\sB$ is finite dimensional, a contradiction. 
Hence $\ker\fe = (0)$, and so $\fe$ is an isomorphism.
\end{proof}

\end{document}